\documentclass[12pt]{article}
\usepackage{graphicx,psfrag,epsfig, color,float}
\usepackage{amssymb,amsmath,amscd,amsthm}
\usepackage{graphicx,psfrag,epsfig}

\usepackage{graphicx}
\usepackage[active]{srcltx}

\newtheorem{theorem}{Theorem}[section]

\def\be{\color{black}}

\date{}

\begin{document}

\date{}
\title{Front propagation for reaction-diffusion equations in composite structures}
\author{
 M. Freidlin\footnote{Dept of Mathematics, University of Maryland,
College Park, MD 20742, mif@math.umd.edu}, L.
Koralov\footnote{Dept of Mathematics, University of Maryland,
College Park, MD 20742, koralov@math.umd.edu}
} \maketitle

\begin{abstract}
We consider asymptotic problems concerning the motion of interface separating the regions of large and
small values of the solution of a reaction-diffusion equation in the media consisting of domains with different characteristics
(composites). Under certain conditions, the motion can be described by the Huygens principle in the appropriate Finsler (e.g., Riemannian) metric.   In general, the motion of the interface has, in a sense, non-local nature. In particular, the interface may move by jumps. 
We are mostly concerned with the nonlinear term that is of KPP type.  The results are based on limit theorems for large deviations. 
\end{abstract}



\section{Introduction.} \label{introd}

Consider a reaction-diffusion equation (RDE)
\begin{equation} \label{rde1itone}
\frac{\partial u^\varepsilon}{\partial t} = \frac{\varepsilon}{2} \sum_{i,j =1}^n a_{ij}(x) \frac{\partial^2 u^\varepsilon}{\partial x_i \partial x_j}   + \frac{1}{\varepsilon} c(x, u^\varepsilon)u^\varepsilon = \varepsilon M u^\varepsilon + \frac{1}{\varepsilon} c(x, u^\varepsilon)u^\varepsilon,~~
t > 0,~x \in \mathbb{R}^n,~~~~ 
\end{equation}
\[
u^\varepsilon(0,x) = g(x) \geq 0.
\]
Here $M$ is an elliptic operator with sufficiently regular coefficients, $\varepsilon > 0$ is a small parameter, and the
nonlinear term is of Kolmogorov-Petrovskii-Piskunov (KPP) type. The latter means that  $c(x, 1) = 0$, 
$c(x, u) < 0$ for $u > 1$, and  $c(x,0) > c(x, u) > 0$ for  $u \in (0,1)$ and $x \in \mathbb{R}^n$.
Assume that $0 \leq g \leq 1$ is  continuous with compact support $G_0$. (We could also allow $g$ to be continuous everywhere except a smooth hypersurface.
In this case, we require that $G_0$ coinsides with the closure of its interior.) 
 We assume that $g$ is not identically equal to zero.  
 We assume that $c$ is Lipschitz continuous in $u$ (uniformly in $x$). 

It was shown in \cite{F1}, \cite{F2}, \cite{F3}  that if $c(x, 0) = \tilde{c}(x) =  \tilde{c}$ is constant, then
$\lim_{\varepsilon \downarrow 0} u^\varepsilon(t,x)$ is equal to zero if $\rho(x, G_0) > t \sqrt{2\tilde{c}}$ and  is equal to one if $\rho(x, G_0) < 
t \sqrt{2\tilde{c}}$,
where $\rho$ is the  Riemannian metric corresponding to the diffusion matrix $a(x) = (a_{ij}(x))$:
\[
\rho(x,y) = \inf_{\substack{\varphi \in C^1([0,1], \mathbb{R}^n) \\ \varphi(0) = x, \varphi(1) = y}} \int_0^1 \sqrt{ (a^{-1}(\varphi(t)) \dot{\varphi}(t), \dot{\varphi}(t))} dt. 
\] 
This result means that when $\varepsilon \ll 1$ the interface between the region where $u^\varepsilon(t,x)$ is close to zero and the region where it is close to one moves according
to the Huygens principle with the constant speed $\sqrt{2\tilde{c}}$ in the metric $\rho$. 

If $\tilde{c}(x)$ is not constant, the position of the interface at 
time $t_2 > t_1$, in general, is not defined by the position of the interface at time $t_1$. Its motion
is, in  a sense, non-local. In particular, it can have jumps (\cite{F1}, \cite{F3}). In general case, the limiting behavior
of $u^\varepsilon(t,x)$ as $\varepsilon \downarrow 0$ 
can be described using the limit theorems for large deviations (see  \cite{FW}).  Let $X^\varepsilon_t$ be the diffusion process on $\mathbb{R}^n$ governed by
the operator $\varepsilon M$:
\begin{equation} \label{sde1}
d X^\varepsilon_t = \sqrt{\varepsilon} \sigma(X^\varepsilon_t) d W_t,~~X^\varepsilon_0 = x,
\end{equation}
where $W_t$ is a Wiener process and $\sigma(x) \sigma^*(x) = a(x)$. The Feynman-Kac formula implies that the solution $u^\varepsilon$ of problem (\ref{rde1itone})
satisfies the following equation
\begin{equation} \label{eqw3}
u^\varepsilon(t,x) = \mathrm{E}_x \left( g( X^\varepsilon_t) \exp\left(\frac{1}{\varepsilon}\int_0^t c(X^\varepsilon_s,
 u^\varepsilon(t-s, X^\varepsilon_s)) ds \right) \right),
\end{equation}
where $\mathrm{E}_x $ means the expected value for trajectories of (\ref{sde1}) with the initial condition $X^\varepsilon_0 = x$. In the case of KPP-type nonlinear
term, (\ref{eqw3}) implies that
\begin{equation} \label{ineq44}
u^\varepsilon(t,x) \leq \mathrm{E}_x \left( g( X^\varepsilon_t) \exp\left(\frac{1}{\varepsilon}\int_0^t \tilde{c}(X^\varepsilon_s) ds \right) \right) = 
\tilde{u}^\varepsilon(t,x).
\end{equation}
Note that the function $\tilde{u}^\varepsilon$ is the solution of the linear problem obtained from (\ref{rde1itone}) when $c(x,u)$ is replaced by $\tilde{c}(x)$. The 
asymptotics of  $\tilde{u}^\varepsilon(t,x)$ in the right hand side of (\ref{ineq44}) can be calculated using large deviation estimates. Namely, if 
$S_{0t}(\varphi)$, $\varphi \in C([0,t], \mathbb{R}^n)$, is the action functional (\cite{FW}) of the family $X^\varepsilon_t$ as $\varepsilon \downarrow 0$
with the normalizing factor $\varepsilon^{-1}$, then
\[
\lim_{\varepsilon \downarrow 0} \varepsilon \ln \tilde{u}^\varepsilon(t,x) = \sup_{\varphi_0 = x, \varphi_t \in G_0} \left( \int_0^t \tilde{c}(\varphi_s) ds
-S_{0t}(\varphi) \right) = \tilde{V}(t,x).
\]
This implies that
\[
\lim_{\varepsilon \downarrow 0} {u}^\varepsilon(t,x)  =  \lim_{\varepsilon \downarrow 0} \tilde{u}^\varepsilon(t,x) = 0 ~~~{\rm if}~~\tilde{V}(t,x) <  0.
\]
Under certain assumptions, one can prove that $\lim_{\varepsilon \downarrow 0} {u}^\varepsilon(t,x)  = 1$ if $\tilde{V}(t,x) >  0$. In this case, the equation
$\tilde{V}(t,x) = 0$ defines the position of the  interface. In particular, if $\tilde{c}(x) = \tilde{c}$ is constant, the position of the interface is
described by the Huygens principle, as above. In the general case, the  position of the interface is 
defined (see \cite{F4}, \cite{FLee}) by the function
\[
{V}(t,x) = \sup_{\varphi_0 = x, \varphi_t \in G_0} \min_{a \in [0,t]} \left( \int_0^a \tilde{c}(\varphi_s) ds
-S_{0a}(\varphi) \right).
\]
If ${V}(t,x) < 0$, then $\lim_{\varepsilon \downarrow 0} {u}^\varepsilon(t,x)  = 0$, while $\lim_{\varepsilon \downarrow 0} {u}^\varepsilon(t,x)  = 1$ if $(t,x)$ belongs
to the interior of the set $\{(t,x): {V}(t,x) = 0 \}$. These results were later re-proved and generalized using classical PDE methods  (see \cite{ES}, \cite{BES}). 

Equation (\ref{eqw3}), together with (\ref{sde1}), is equivalent to (\ref{rde1itone}). It describes the interplay between the transport of particles (in our case the diffusion of particles) and the law of multiplication/annihilation of particles. Note that, instead of the diffusion transport defined by (\ref{sde1}), one could
consider other types of stochastic motion, as long as the action functional for the family is known and a certain Markov property is satisfied. One could also consider a non-local non-linear term (compare with \cite{AFV}).


In this paper, we will study interface propagation  for reaction-diffusion equations in 
composite structures. By a composite structure we mean a domain that is a union of two or more regions with significantly different
properties of the media (coefficients of the equation). In the case of  layered structures that are space-homogeneous (in each of the layers), it turns out that 
the interface motion can also be described by the Huygens principle. However, the speed of the motion is constant if it is calculated with respect to 
an appropriate Finsler metric, rather than a Riemannian metric. We derive the expression for this metric in three qualitatively different cases, depending 
on the magnitude of the underlying diffusion across the layers.

In contrast to the case of a single layer, now the propagation of the
interface is not described by the Huygens principle and may be non-local, even if the nonlinear term does not vary within each of the layers. 
The main difference between the case of the single layer and the one with several layers is that now the propagation of the interface is determined not only
by the large deviations of the underlying diffusion along the layer, 
but by the interplay between the deviation from the stationary destribution between the layers  and the 
large deviations for the diffusion in each of the layers. A similar, in a sense, phenomenon was studied in \cite{FGMP}.
\vskip -15pt
\begin{figure}[htbp]
    \centerline{\includegraphics[height=3.8in, width= 6.5in,angle=0]{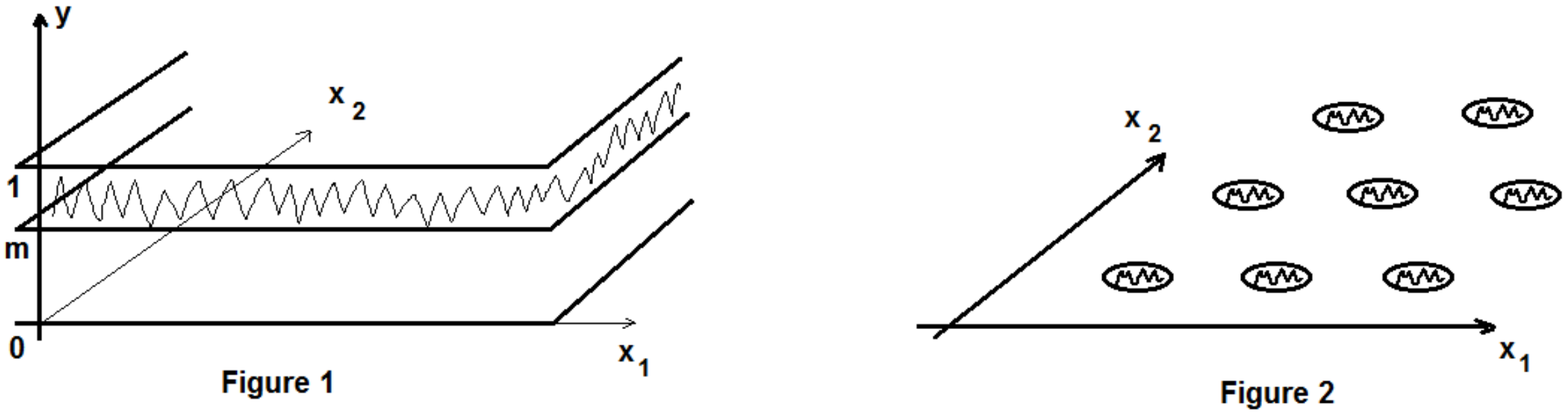}}
\vskip -150pt
    \label{secondkind}
\end{figure}

Examples of composite structures are given in   Figures 1  and 2. 
The composite in  Figure 1  consists of two layers with different properties.   Figure 2  shows periodic 
inclusions in  a homogeneous medium. 
First, let us consider the layered structure shown in  Figure~1. 

The reaction-diffusion equation in a structure with two layers has the form
\begin{equation} \label{rde1it}
\frac{\partial u^\varepsilon}{\partial t} = \frac{\varepsilon}{2} \sum_{i,j =1}^n a_{ij}(x,y) \frac{\partial^2 u^\varepsilon}{\partial x_i \partial x_j}  + 
\frac{\varepsilon^{-\beta}}{2} \alpha(x,y)
\frac{\partial^2 u^\varepsilon}{\partial y^2} + \frac{1}{\varepsilon} c(x, y, u^\varepsilon)u^\varepsilon,~
t > 0,~x \in \mathbb{R}^n,~y \in I_1 \cup I_2,
\end{equation}
\[
\frac{\partial u^\varepsilon}{\partial y}|_{y=0,1} = 0,~~~u^\varepsilon(0,x,y) = g(x),
\]
where $0 \leq g \leq 1$ is  continuous with compact support $G_0$ and is not identically equal to zero, 
 $I_1 = (0,m)$, and $I_2 = (m,1)$. It is assumed that $g$ is not identically equal to zero.
To account for different layers, we  assume that
\[
a_{ij}(x,y) = \begin{cases} a_{ij}^1(x), ~ y \in I_1\\ a_{ij}^2(x), ~  y \in I_2, \end{cases}~\alpha(x,y) = \begin{cases} \alpha^1(x), ~  y \in I_1 \\ 
\alpha^2(x), ~ y \in I_2. \end{cases}
\]
It is assumed that $a^1$,$a^2$ are uniformly bounded and uniformly positive-definite matrices and that $\alpha^1$, $\alpha^2$ are uniformly bounded and uniformly positive. 
The nonlinear term in (\ref{rde1it}) also depends on the layer: we assume that
\[
c(x, y,u) = 
\begin{cases} c^1(x,u), ~  y \in I_1, \\ c^2(x,u), ~ y \in I_2. \end{cases}
\]
It is assumed that $a_{ij}^k, \alpha^k \in C^2(\mathbb{R}^n)$, while $c^k(x,u)$ is Lipschitz continuous,  $k =1,2$. 
Observe that the coefficients in (\ref{rde1it}) may be discontinuous for $y = m$, and the equation is satisfied only
when $y \neq m$. Thus, in order for the uniqueness of the solution to hold, one should add a `gluing condition' on the plane $y = m$. 
To do this rigorously, it is best to relate $u^\varepsilon$ to itself using the Feynman-Kac formula, similarly to (\ref{eqw3}), and then use this
as the definition of the solution of (\ref{rde1it}) (compare with \cite{F3}).  Namely, let
\[
M^\varepsilon u(x,y) = \frac{\varepsilon}{2} \sum_{i,j =1}^n a_{ij}(x,y) \frac{\partial^2 u}{\partial x_i \partial x_j}  + 
\frac{\varepsilon^{-\beta}}{2} \alpha(x,y)
\frac{\partial^2 u}{\partial y^2}.
\]
The domain of $M^\varepsilon$ contains functions $u \in C(\mathbb{R}^n \times [0,1]) \bigcap C^2(\mathbb{R}^n \times (I_1 \bigcup I_2))$, whose
first derivative in $y$ belongs to $C(\mathbb{R}^n \times [0,1])$, which satisfy $\frac{\partial u}{\partial y}|_{y=0,1} = 0$, and are such that 
$M^\varepsilon u$  (understood as the differential operator on $\mathbb{R}^n \times (I_1 \bigcup I_2)$ applied to $u$ and extended to $\mathbb{R}^n \times [0,1]$)
belongs to $C(\mathbb{R}^n \times [0,1])$. The closure of the operator $M^\varepsilon$ with the domain specified above serves as the generator for the Markov family 
$(X^{\varepsilon}_t, Y^{\varepsilon}_t)$ on $\mathbb{R}^n \times [0,1]$ (where we dropped the  dependence on the initial point $(x,y)$ from the notation). 
This diffusion process is the  limit of processes with continuous diffusion coefficients approximating  the
diffusion coefficients $a_{ij}(x,y)$. 

The pair of processes $(X^{\varepsilon}_t, Y^{\varepsilon}_t)$ starting at $(x,y)$ is  the solution of the system of
stochastic differential equations 
\[
d X^{\varepsilon}_t = \sqrt{\varepsilon} A(X^{\varepsilon}_t, Y^{\varepsilon}_t) dW_t,~~~X^{\varepsilon}_0 = x,
\]
\[
d Y^{\varepsilon}_t = \frac{1}{\sqrt{\varepsilon^\beta}} \sigma (X^{\varepsilon}_t, Y^{\varepsilon}_t) dV_t,~~~Y^{\varepsilon}_0 = y,
\]
where $A$ is the positive-definite symmetric square root of the matrix $a$, $\sigma = \sqrt{\alpha}$,  $W_t$ is an $n$-dimensional Brownian motion, and $V_t$ is a one-dimensional Brownian motion  independent of $W_t$. 
The process $Y^{\varepsilon}_t$ is reflected at the end points of the segment and
satisfies a gluing condition at $y = m$. We define the solution of  (\ref{rde1it}) as the bounded continuous function on $[0,\infty) \times \mathbb{R}^n \times [0,1]$
that sastisfies
\[
u^\varepsilon(t,x,y) =
\]
\[
 \mathrm{E}_{(x,y)} \left( g(X^{\varepsilon}_{t}) \exp(\varepsilon^{-1} \int_0^{t} c(Y^{\varepsilon}_s, 
u^\varepsilon(t -s, X^{\varepsilon}_s, Y^{\varepsilon}_s) ds) \right)
\]
for each $t,x,y$. Using the Lipschitz continuity of $c$ in the second argument, it is easy to show that
such a function $u^\varepsilon$ exists and is unique  for each $\varepsilon > 0$. 


We are  mostly  interested in the case when the nonlinearity is of Kolmogorov-Petrovskii-Piskunov (KPP) type.  Namely, we assume that $c^k$, $k =1,2$, are  uniformly  Lipschitz continuous in $u$,  $c^k(x, 1) = 0$, 
$c^k(x, u) < 0$ for $u > 1$, and  $c^k(x,0) > c^k(x, u) > 0$ for $x \in \mathbb{R}^n$  and $u \in (0,1)$.  

The asymptotics of $u^\varepsilon$ as $\varepsilon \downarrow 0$ for various values of the parameter $\beta$ will be studied in this paper. In Section~\ref{lin1aa},
we consider the situation when $a^k$ and $c^k$, $k =1,2$, do not depend on the $x$-variable. In this case, the metric governing the interface propagation is
translation-invariant - it is given by a norm of the difference between the points. 
Three different cases are distinguished, depending on whether $\beta = 1$, $\beta > 1$, or $\beta < 1$. In Section~\ref{lin1bb}, we discuss the situation when $a^k$,
$\alpha^k$,  
and $c^k$ are allowed to depend on $x$. In both Section~\ref{lin1aa} and Section~\ref{lin1bb}, we use the large deviation 
principle for the joint distribution of the trajectory of the underlying diffusion in the $x$-space and the occupation measure for the diffusion in the $y$-space.
In the case of $x$-dependent coefficients, the large deviation principle is more complicated.  

%

\section{The case of $x$-independent coefficients.} \label{lin1aa}
\subsection{Asymptotics of solutions to linear equations.}  \label{lint}

In this section, we consider the linear version of the Cauchy problem (\ref{rde1it}), i.e., we assume that $c(x,y,u) = c(x,y)$. 
The coefficients
$a^k, \alpha^k, c^k$, $k =1,2$, do not depend on~$x$. Thus $u^\varepsilon$ satisfies
\begin{equation} \label{rde1itlin}
\frac{\partial u^\varepsilon}{\partial t} = \frac{\varepsilon}{2} \sum_{i,j =1}^n a_{ij}(y) \frac{\partial^2 u^\varepsilon}{\partial x_i \partial x_j}  + 
\frac{\varepsilon^{-\beta}}{2} \alpha(y)
\frac{\partial^2 u^\varepsilon}{\partial y^2} + \frac{1}{\varepsilon} c(y)u^\varepsilon,~~
t > 0,~x \in \mathbb{R}^n,~y \in (0,1) \setminus \{m\}. 
\end{equation}
\[
\frac{\partial u^\varepsilon}{\partial y}|_{y=0,1} = 0,~~~u^\varepsilon(0,x,y) = g(x).
\]
 We will show that there is a function $ \lambda(t, x)$, continuous on $(0,\infty) \times \mathbb{R}^n$, such that
$\varepsilon \ln {u}^\varepsilon (t,x,y) \rightarrow  \sup_{x' \in G_0} \lambda(t, x - x')$. 
 The expressions for $\lambda(t,x)$ are different, depending on whether $\beta = 1$, $\beta > 1$, or $-1 < \beta < 1$. (If $\beta = - 1$, there is
no need to distinguish between the $x$ and $y$ variables, and the answer follows from \cite{F3}. If $\beta < - 1$, then, in order to find the asymptotics 
of ${u}^\varepsilon (t,x,y)$ with $y \neq m$,  the equation can be viewed in the $(t,x)$ space, with the diffusion in the $y$ variable ignored, and the $y$ variable in the coefficients treated as a parameter.) The function $\lambda$ is the multi-layer analogue of $\tilde{V}$ defined in the Introduction. 

First consider the case when $\beta = 1$.  Let 
\[
Lu (y) = \frac{1}{2} \alpha(y)
 u''(y)
\]
be the operator on $C([0,1])$ with the domain $D(L)$ that  consists of functions satisfying  
\[
u \in C^1([0,1]),~~~ \alpha u'' \in C([0,1]),~~~u'(0) = u'(1) = 0.
\]
Let $Y^{\varepsilon}_t$ be the process s with values on $[0,1]$, whose generator is $ \varepsilon^{-\beta} L$.
Thus, if the initial value of the process  $Y^{\varepsilon}_t$ is $y$, the process formally satisfies
\[
d Y^{\varepsilon}_t = \frac{1}{\sqrt{\varepsilon^\beta}} \sigma (Y^{\varepsilon}_t) dV_t,~~~Y^{\varepsilon}_0 = y,
\]
where $\sigma = \sqrt{\alpha}$ and $V_t$ is a one-dimensional Brownian motion. ($Y^{\varepsilon}_t$ is reflected at the end points of the segment and
satisfies a gluing condition at $y = m$.)

Given initial values $X^{\varepsilon}_0 = x$ and $Y^\varepsilon_0 = y$, define 
\[
X^{\varepsilon}_t = x + \sqrt{\varepsilon}  \int_0^t A(Y^{\varepsilon}_s) dW_s,
\]
where $A$ is the positive-definite symmetric square root of the matrix $a = (a_{ij})$ and $W_t$ is an $n$-dimensional Brownian motion independent of $V_t$. Note
that $X^{\varepsilon}_t$ also depends on $\beta$, although this is not reflected in the notation. 

We will repeatedly make use of the following simple observation (compare with \cite{FW}, Ch. 3). Let $\Lambda^\varepsilon_z$ be a family of probability measures on $(M, \mathcal{B}(M))$, where
$(M,d)$ is a metric space, $\varepsilon > 0$ is a small parameter, and $z$ is an additional parameter  
(for example, $\Lambda^\varepsilon_z$ may be the measures induced by processes that start at an initial point $z$). 
Suppose that $S_z$ is  the action functional for 
$\Lambda^\varepsilon_z$ with  normalizing coefficient $\varepsilon^{-1}$, uniformly in $z$. 
Then for continuous functions $0 \leq \varphi \leq C$ and $\psi \leq C$ on~$M$,
\begin{equation} \label{actf}
\lim_{\varepsilon \downarrow 0} \left( \varepsilon \ln \int_M \varphi(x) \exp({ \frac{\psi(x)}{\varepsilon} }) d \Lambda^\varepsilon_z \right) 
= \sup_{x \in {\rm supp}(\varphi)} (\psi(x)-S_z(x)),
\end{equation}
uniformly in $z$. 

If $\psi$ is not continuous, we can still estimate the left-hand side of (\ref{actf}) from above. Namely, for $\eta > 0$ define $S^\eta_z(x) = \inf_{y: d(y,x) \leq \eta} S_z(y)$. Then it is not difficult to see that
\begin{equation} \label{actf2}
\lim_{\varepsilon \downarrow 0} \left( \varepsilon \ln \int_M \varphi(x) \exp({ \frac{\psi(x)}{\varepsilon} }) d \Lambda^\varepsilon_z \right) 
\leq \sup_{x \in {\rm supp}(\varphi)} (\psi(x)-S^\eta_z(x)),
\end{equation}
uniformly in $z$.

For $f \in C([0,1])$, let $H(f)$ be the top eigenvalue 
of the operator $L_f u = Lu + fu$. Let $\mathcal{M}_{[0,1]}$ be the space of probability measures on $([0,1], \mathcal{B}([0,1]))$.
 Let $\mathcal{M}'_{[0,1]} = \{ \mu \in \mathcal{M}_{[0,1]}: \mu(\{m\}) = 0 \}$. Let $\mu^{\varepsilon}_{t,y}$ be the normalized occupation measure on 
$([0,1], \mathcal{B}([0,1]))$ of the process 
$Y^{\varepsilon}$ (with $Y^{\varepsilon}_0 = y$) on the interval $[0, t]$, i.e., $\mu^{\varepsilon}_{t,y}(B) = \int_0^t \chi_B(Y^{\varepsilon}_s) ds /t$, 
$B \in \mathcal{B}([0,1])$.

For $\mu \in \mathcal{M}_{[0,1]}$, define 
\[
I(\mu) = \sup_{f \in C([0,1])} (\int_0^1 f d \mu - H(f)). 
\]
Then $t I $ is the action functional for $\mu^{\varepsilon}_{t,y}$, uniformly in $(t, y) \in [a,b] \times [0,1]$ if $0 < a < b$ (see \cite{Gartner}, \cite{FW} (Ch. 10)). 
Let 
\[
J = \{p = (p_1, p_2): p_1 + p_2 =1,~p_1, p_2 \geq 0 \}.
\]
This space is endowed with the metric  $d_J((p'_1, p'_2), (p''_1, p''_2)) = |p'_1 - p''_1|$.  
For $p \in J$ and $\mu \in \mathcal{M}'_{[0,1]}$, define $p_\mu = (\mu(I_1), \mu(I_2))$ and
\begin{equation} \label{defes}
S(p) = \inf_{\mu: p_\mu = (p_1, p_2)} I(\mu).
\end{equation}
Thus $t S$ is the action functional,  uniformly in $(t, y) \in [a,b] \times [0,1]$,  for the family of measures on $J$ induced by the random vectors 
$(\mu^{\varepsilon}_{t,y}(I_1), \mu^{\varepsilon}_{t,y}(I_2))$. 
Such measures (which also depend on $\beta$) will be denoted by $\Lambda^\varepsilon_{t,y}$, i.e., 
\[
\Lambda^\varepsilon_{t,y} (A) = \mathrm{P} (p_{\mu^{\varepsilon}_{t,y}} \in A),~~~A \in \mathcal{B}(J). 
\] 

 In order to derive the asymptotics of  ${u}^\varepsilon (t,x,y)$, we will
show that the main contribution to the expectation in the Feynman-Kac formula comes from the event where the trajectories of the underlying diffusion spend an asymptotically non-random proportion of time  $p_1$ in the region where $y \in I_1$,
and an asymptotically non-random proportion of time  $p_2$ in the region where $y \in I_2$. Assuming that $p_1$ and $p_2$ are known, we will derive 
the expression for the contribution to the expectation in the Feynman-Kac formula, and then maximize the expression under the condition that $p_1+ p_2 = 1$.

 Let $a^1 = (a^1_{ij}), a^2 = (a^2_{ij})$. 
For $v \in \mathbb{R}^n$, define
\[
R(p, v) = \frac{1}{2} ((p_1 a^1 + p_2 a^2)^{-1} v, v),
\]
\[
T(p) =  p_1 c^1 + p_2 c^2. 
\]
Now we can write the expression for $\lambda(t,x)$ in the case when $\beta = 1$,
\begin{equation} \label{mo}
\lambda(t,x) = \sup_{p} (t (T(p) - S(p)  - {R(p, \frac{x}{t})})). 
\end{equation}

Next consider the case $\beta > 1$.  The difference from the case with $\beta = 1$ is that now the values of $p_1$ and
$p_2$ are prescribed. Namely, let $\pi$ be the invariant measure
for the process $Y^{y,\varepsilon}_t$ (the invariant measure doesn't depend on $\varepsilon$ or $\beta$). 
The expression for $\lambda(t,x)$ in the case when $\beta > 1$ is
\begin{equation} \label{ltmo}
\lambda(t,x) = t( T(p_\pi)  - {R(p_\pi, \frac{x}{t})}).
\end{equation}

Finally, consider $-1 < \beta < 1$.  In this case, we again have minimization in $p$, but the term $S(p)$ is not present. Namely, define
\begin{equation} \label{gtmo}
\lambda(t,x) = \sup_{p} (t (T(p) - {R(p, \frac{x}{t})})). 
\end{equation}

\begin{theorem} \label{linle}
Under the above assumptions, 
\begin{equation} \label{limyy}
\lim_{\varepsilon \downarrow 0} \varepsilon \ln ({u}^\varepsilon(t,x,y)) =  \sup_{x' \in G_0} \lambda(t, x - x')
\end{equation}
uniformly on every compact $K \subset (0,\infty) \times \mathbb{R}^n \times [0,1]$, where $\lambda$ is given by (\ref{mo}) if $\beta = 1$, by 
(\ref{ltmo}) if $\beta > 1$, and (\ref{gtmo}) if $-1 < \beta < 1$. 
\end{theorem}
\proof  
Fix $t > 0$. Let $M^{x,p, \varepsilon}_t$ be the measure on $\mathcal{C} = C([0,t], \mathbb{R}^n)$ induced by the process~$X^{\varepsilon}_t$ conditioned
on $\Lambda^{\varepsilon}_{t,y}(\{p\}) = 1$ (obesrve that there is no dependence on $y$ or $\beta$ in $M^{x,p, \varepsilon}_t$, as follows from the definition of the 
process $X^{\varepsilon}_t$).

By the Feynman-Kac formula, 
\begin{equation} \label{fkfo1}
u^\varepsilon(t,x,y) = \mathrm{E}_{(x,y)} \left( g(X^{\varepsilon}_t) \exp(\varepsilon^{-1} \int_0^t c(Y^{\varepsilon}_s) ds) \right)= 
\end{equation}
\[
\int_J \exp (\varepsilon^{-1} t (c^1 p_1 + c^2 p_2) ) \int_{\mathcal{C}} g(\varphi_t)    d M^{x, p, \varepsilon}_t(\varphi) d {\Lambda}^{\varepsilon}_{t,y}(p). 
\]
For a compact $\bar{K} \subset \mathbb{R}^n$, the action functional for $M^{x, p, \varepsilon}_t$ is given, uniformly in $(x,p) \in \bar{K} \times J$,  by 
$\int_0^t R(p, \varphi'(s)) ds$ when $\varphi(0) = x$ (and is equal to $-\infty$ otherwise). Therefore, by (\ref{actf}),
\[
\lim_{\varepsilon \downarrow 0} \left( \varepsilon \ln \int_{\mathcal{C}} g(\varphi_t)    d M^{x, p, \varepsilon}_t(\varphi) \right) = - \inf_{\varphi: \varphi(0) = x,
\varphi(t) \in G_0} \int_0^t R(p, \varphi'(s)) ds = 
-\inf_{x' \in G_0} t R(p, \frac{x - x'}{t}),
\]
uniformly in $(x,p) \in \bar{K} \times J$. Substituting this in (\ref{fkfo1}), we get
\begin{equation} \label{nkk}
\lim_{\varepsilon \downarrow 0} \varepsilon \ln ({u}^\varepsilon(t,x,y)) =  \lim_{\varepsilon \downarrow 0} \varepsilon \ln 
\int_J \exp \left(\varepsilon^{-1} t (T(p)  -\inf_{x' \in G_0} R(p, \frac{x - x'}{t})) \right) d {\Lambda}^{\varepsilon}_{t,y}(p).
\end{equation}
When $\beta = 1$, we use (\ref{actf}) and the fact that $tS$ is the action functional for the family ${\Lambda}^{\varepsilon}_{t,y}$ in order to obtain
\[
\lim_{\varepsilon \downarrow 0} \varepsilon \ln ({u}^\varepsilon(t,x,y)) = \sup_{x' \in G_0} \sup_{p} (t (T(p) - S(p)  - {R(p, \frac{x-x'}{t})})),
\]
uniformly in $(x, y) \in \bar{K} \times [0,1]$. Next, consider the case when $\beta > 1$.  If $U \subseteq J$ is an open neighborhood of $p_\pi$, then, for each $C > 0$,
\[
\Lambda^{\varepsilon}_{t,y}(U) \geq 1 - \exp(- \varepsilon^{-1} C)
\]
for all sufficiently small $\varepsilon$. Therefore, the main contribution to the integral in (\ref{nkk})  comes from an arbitrarily small
neighborhood of $p_\pi$, which implies that
\[
\lim_{\varepsilon \downarrow 0} \varepsilon \ln ({u}^\varepsilon(t,x,y)) = \sup_{x' \in G_0} t( T(p_\pi)  - {R(p_\pi, \frac{x-x'}{t})}),
\]
uniformly in $(x, y) \in \bar{K} \times [0,1]$.
Finally, if $-1 < \beta < 1$, then for each nonempty open set $U \subseteq J$ and each $c > 0$ we have
\[
\Lambda^{\varepsilon}_{t,y}(U) \geq \exp(- \varepsilon^{-1} c)
\]
for all sufficiently small $\varepsilon$. Therefore, 
\[
\lim_{\varepsilon \downarrow 0} \varepsilon \ln ({u}^\varepsilon(t,x,y)) = \sup_{x' \in G_0} \sup_{p} (t (T(p) - {R(p, \frac{x-x'}{t})})),
\]
uniformly in $(x, y) \in \bar{K} \times [0,1]$. We have thus justified (\ref{gtmo}) in all the three cases for fixed~$t > 0$. Let us now show that the convergence is
uniform on $K \subset (0,\infty) \times \mathbb{R}^n \times [0,1]$. From the Feynman-Kac formula  it follows that for $\delta < t$
\[
\mathrm{E}_{(x,y)} \left( u^\varepsilon(t-\delta, X^{\varepsilon}_\delta, Y^{\varepsilon}_\delta) \right) \leq u^\varepsilon(t,x,y) \leq \exp(\frac{\delta \max(c^1, c^2)}{\varepsilon} )\mathrm{E}_{(x,y)} \left( u^\varepsilon(t-\delta, X^{\varepsilon}_\delta, Y^{\varepsilon}_\delta) \right).
\]
Considering the contribution to the expectation from the events $\|X^{\varepsilon}_\delta - x\| \leq \eta$ and $\|X^{\varepsilon}_\delta - x\| > \eta$ and using
the large deviations estimates on the process $X^{\varepsilon}_t$, we see that for each $\eta > 0$ and $\alpha > 0$ there exist $\delta_0 > 0, \varepsilon_0 > 0$ such that
\[
\frac{1}{2} \inf_{x': \|x' -x\| \leq \eta} \inf_{y \in [0,1]}  u^\varepsilon(t-\delta, x', y)  \leq
\]
\[
 u^\varepsilon(t,x,y) \leq 
\]
\[
\exp(\frac{\delta}{\varepsilon} )\sup_{x': \|x' -x\| \leq \eta} \sup_{y \in [0,1]}  u^\varepsilon(t-\delta, x', y) + 
\exp(\frac{t \max(c^1, c^2)-\alpha}{\varepsilon}),
\]
when $\delta < \delta_0$, $\varepsilon < \varepsilon_0$. Together with the convergence in (\ref{gtmo})  for fixed~$t > 0$ and the continuity of the right hand
side of (\ref{gtmo}), this is enough to conclude that  the convergence is
uniform on $K \subset (0,\infty) \times \mathbb{R}^n \times [0,1]$.

\qed
\\
\\
{\bf Remark.}  In the proof of Theorem~\ref{linle} we saw that for each $r, \delta > 0$
\[
\varepsilon \ln \mathrm{E}_{(x,y)} \left( g(X^{\varepsilon}_t) \exp(\varepsilon^{-1} \int_0^t c(Y^{\varepsilon}_s) ds),~X^{\varepsilon}_t 
\in B_r(x_0) \right) \geq
\lambda(t, x - x_0) - \delta
\]
for all sufficiently small $\varepsilon$, where $x_0 \in {\rm Int}(G_0)$. The same argument gives the bound if we restrict the expectation to the event that 
$X^{\varepsilon}_t$ closely follows the segment connnecting
$x$ to $x_0$. More precisely, let $\varphi: [0,t] \rightarrow \mathbb{R}^n$ be the linearly parametrized segment with $\varphi(0) = x$, $\varphi(t) = x_0$. Then
\begin{equation} \label{lpp}
\varepsilon \ln \mathrm{E}_{(x,y)} \left( g(X^{\varepsilon}_t) \exp(\varepsilon^{-1} \int_0^t c(Y^{\varepsilon}_s) ds),~
\sup_{s \in [0,t]}\|X^{\varepsilon}_s - \varphi(s)\| \leq 
\delta \right) \geq
\lambda(t, x - x_0) - \delta
\end{equation}
for all sufficiently small $\varepsilon$, uniformly on every compact $K \subset (0,\infty) \times \mathbb{R}^n \times [0,1]$.
\\
\\
{\bf Remark.} For $t_1, t_2 \geq 0$, we have 
\begin{equation} \label{conv1}
\lambda(t_1+t_2, x_1+x_2) \geq  \lambda(t_1, x_1) + \lambda(t_2, x_2).
\end{equation}
Indeed, suppose that $G_0 = B_{r_0}(0)$ (the ball of radius $r_0$ around the origin). 
By Theorem~\ref{linle}, for each $\delta > 0$ there is $r> 0$ such that, for 
all sufficiently small~$\varepsilon$, 
\[
{u}^\varepsilon(t_1, x', y') \geq \exp((\lambda (t_1, x_1) - \delta)/\varepsilon)
\]
when $x' \in B_r(x_1)$. Let $\tilde{g}$ be a continuous function taking values in $[0,1]$ that is equal to one
on $B_{r/2}(x_1)$ and equal to zero outside $B_r(x_1)$. Applying Theorem~\ref{linle} again, this time on the interval $[0, t_2]$, with initial function $\tilde{g}$,
and using the semigroup property of solutions to the linear equation, we obtain
\[
{u}^\varepsilon(t_1 + t_2,x_1+x_2,y) \geq \exp((\lambda (t_1, x_1) - \delta)/\varepsilon) \exp((\lambda (t_2, x_2) - \delta)/\varepsilon), 
\]
and therefore
\[
\varepsilon \ln {u}^\varepsilon(t_1 + t_2,x_1+x_2,y) \geq \lambda(t_1, x_1) + \lambda(t_2, x_2) - 2 \delta.
\]
The left hand side can be made arbitrarily close to $\lambda(t_1+t_2, x_1+x_2)$ by selecting a sufficiently small $r_0$ and a sufficiently small $\varepsilon$. 
Thus, since $\delta > 0$ was arbitrary, we obtain~(\ref{conv1}). 

\subsection{Asymptotics of solutions to reaction-diffusion equations.}  \label{rde1}

In this section we consider the Cauchy problem for the reaction-diffusion equation (\ref{rde1it}). It is assumed that $a^k,  c^k$, $k =1,2$, do not depend on~$x$.
Thus
\[
c(y,u) = 
\begin{cases} c^1(u), ~  y \in (0,m) \\ c^2(u), ~ y \in (m,1). \end{cases}
\]
Let $\tilde{c}^1 = c^1(0)$, $\tilde{c}^2 = c^2(0)$. 
Consider the linear problem (\ref{rde1itlin})
with $c^1$, $c^2$ replaced by $\tilde{c}^1$, $\tilde{c}^2$. Let $\lambda(t,x)$ be
given by (\ref{mo}) if $\beta = 1$, by 
(\ref{ltmo}) if $\beta > 1$, and (\ref{gtmo}) if $-1 < \beta < 1$. 

  Define the norm $\|x\|$ via the condition
\[
\lambda(\| x \|, x) = 0.  
\]
From the definition of $\lambda$, in each of the cases it follows that $\lambda(|a| t, a x) = |a| \lambda(t,x)$, and therefore
$\|ax\| = |a| \|x\|$. The triangle inequality follows from (\ref{conv1}), and so $\|\cdot\|$ is indeed a norm. Let $d(x_1, x_2) =\|x_1 - x_2\|$. 
Define
\[
G_t = \{x \in \mathbb{R}^n: d(x, G_0) \leq t\}. 
\]
Note that the growth of $G_t$ is described by the Huygens principle in the (translation-invariant) metric $d$. 
\begin{theorem} \label{nnln} If $u^\varepsilon(t,x,y)$ is the solution of (\ref{rde1it}) and $c$ is of KPP type, then, for each $t > 0$, 
\[
\lim_{\varepsilon \downarrow 0} u^\varepsilon(t,x,y) = 0
\]
uniformly on every compact $K \subset (\mathbb{R}^n \setminus G_t) \times [0,1]$, and
\[
\lim_{\varepsilon \downarrow 0} u^\varepsilon(t,x,y) = 1
\]
uniformly on every compact $K \subset {\rm Int}(G_t) \times [0,1]$.

\end{theorem}
\proof Let $\tilde{u}$ be the solution of the linear problem (\ref{rde1itlin})
with $c^1$, $c^2$ replaced by $\tilde{c}^1$, $\tilde{c}^2$. Since $c^k(u) \leq \tilde{c}^k$ for $0 \leq u \leq 1$  (the nonlinearity
is of KPP type), it is clear that $u \leq \tilde{u}$.  By Theorem~\ref{linle}, $\lim_{\varepsilon \downarrow 0} \tilde{u}^\varepsilon(t,x,y) = 0$
uniformly on every compact $K \subset (\mathbb{R}^n \setminus G_t) \times [0,1]$, and therefore $\lim_{\varepsilon \downarrow 0} {u}^\varepsilon(t,x,y) = 0$
uniformly on $K$. 

Now consider a compact $K$ such that $K \subset {\rm Int}(G_t) \times [0,1]$. 
Let $\eta > 0$ and $(x_0, y_0) \in K$. Assume that $x_0 \notin G_0$. Let $t_0 = d(x_0, G_0) < t$.   
Let us show that there is $\delta > 0$ such that 
\begin{equation} \label{ebel}
{u}^\varepsilon(t_0,x,y) \geq  \exp(-\varepsilon^{-1} \eta)
\end{equation}
 for all sufficiently small~$\varepsilon$  when $\|x - x_0\|\leq \delta$.  Let $\tilde{\lambda} (t, x) = \sup_{x' \in G_0} \lambda(t, x - x')$. 
Given $\delta_1 > 0$, we can choose $x_1 \in 
{\rm Int}(G_0)$ and $\delta > 0$ in such a way that for each $x$  we have $\tilde{\lambda}(t_0 -\delta_1 - s, \varphi(s)) < 0$ for $s \leq t_0 -2\delta_1$, where
$\varphi$ is  the linearly parametrized segment $\varphi: [0,t_0 - \delta_1] \rightarrow \mathbb{R}^n$, $\varphi(0) = x$, $\varphi(t_0) = x_1$,. Taking, if necessary, a smaller value of $\delta$, we can make sure
that  $\tilde{\lambda}(t_0 -\delta_1 - s, \psi(s)) < -\delta$ for $s \leq t_0 -2\delta_1$ whenever $\psi: [0,t_0 - \delta_1] \rightarrow \mathbb{R}^n$ is such that $\|\varphi(s) -
\psi(s) \| \leq \delta $ for all $s$. Let $\hat{\varphi}: [0,t_0] \rightarrow \mathbb{R}^n$ be defined via 
\[
\hat{\varphi}(s) = 
\begin{cases} x, ~  s \in [0,\delta_1] \\ \varphi(s -\delta_1), ~ s \in [\delta_1, t_0]. \end{cases}
\]
By the Feynman-Kac formula (which defines the solution),
\[
u^\varepsilon(t_0,x,y) =
\]
\begin{equation} \label{lfo}
 \mathrm{E}_{(x,y)} \left( g(X^{\varepsilon}_{t_0}) \exp(\varepsilon^{-1} \int_0^{t_0} c(Y^{\varepsilon}_s, 
u^\varepsilon(t_0 -s, X^{\varepsilon}_s, Y^{\varepsilon}_s) ds) \right) \geq
\end{equation}
\[
\mathrm{E}_{(x,y)} \left( g(X^{\varepsilon}_{t_0}) \exp(\varepsilon^{-1} \int_0^{t_0} c(Y^{\varepsilon}_s, 
u^\varepsilon(t_0 -s, X^{\varepsilon}_s, Y^{\varepsilon}_s) ds),\sup_{s \in [0,t_0]} \|X^{\varepsilon}_s - \hat{\varphi}(s)\| \leq 
\delta \right). 
\]
  Observe that $\tilde{u}^\varepsilon(t_0 -s, x, y) \rightarrow 0$
uniformly for $(s, x)$ such that $\delta_1 \leq s \leq t_0 -  \delta_1$, $\|x- \hat{\varphi}(s)\| \leq \delta$, since 
$\varepsilon \ln \tilde{u}^\varepsilon(t_0 -s, x, y) \rightarrow \tilde{\lambda}(t_0 -s , x) < -\delta$. 
Since $u^\varepsilon \leq
\tilde{u}^\varepsilon$, the right hand side of (\ref{lfo})
can be estimated from below, for all sufficiently small $\varepsilon$,  by 
\[
\mathrm{E}_{(x,y)} \left( g(X^{\varepsilon}_{t_0}) \exp(\varepsilon^{-1} \int_{\delta_1}^{t_0 -\delta_1} 
(\tilde{c}(Y^{\varepsilon}_s) -\delta) ds),~\sup_{s \in [0,t_0]} \|X^{\varepsilon}_s - \hat{\varphi}(s)\|\right). 
\]
Conditioning on the value of the process at time $\delta_1$, we estimate the value of this expression, from below, by the product  $R_1 \times R_2 \times R_3$, where
\[
R_1 =\exp(-\varepsilon^{-1}(\delta t_0 + \delta_1 \max(\tilde{c}^1 + \tilde{c}^2) )),
\]
\[
R_2 = \mathrm{P}_{(x,y)} (\sup_{s \in [0,\delta_1]} \|X^{\varepsilon}_s - x\| \leq 
\frac{\delta}{2}),
\]
\[
R_3 =  \inf_{x', y': \|x' - x\| \leq \frac{\delta}{2}}  \mathrm{E}_{(x',y')} \left( g(X^{\varepsilon}_{t_0 - \delta_1}) 
\exp(\varepsilon^{-1} \int_{0}^{t_0 -\delta_1} 
\tilde{c}(Y^{\varepsilon}_s) ds),\sup_{s \in [0,t_0 -\delta_1]} \|X^{\varepsilon}_s - {\varphi}(s)\| \leq 
\delta \right).
\]

It follows from (\ref{lpp}) that for all sufficiently small $\delta_1$ and $\delta$ (which may depend on $\delta_1$),
\[
\varepsilon \ln R_3  \geq  \inf_{x': \|x' - x\| \leq \frac{\delta}{2}} \tilde{\lambda}(t_0 -\delta_1, x' - x_1) - \frac{\eta}{4}  \geq - \frac{\eta}{2},
\]
provided that $\varepsilon$ is sufficiently small. Also, for all sufficiently small $\delta_1$ and $\delta$ we have
\[
R_1, R_2 \geq \exp(- \frac{ \varepsilon^{-1} \eta}{4}),
\]
for all sufficiently small~$\varepsilon$. Thus $u^\varepsilon(t_0,x,y)$
 can
be made larger than $\exp(-\varepsilon^{-1} \eta)$ for all sufficiently small~$\varepsilon$.

Now suppose that  $x_0 \in G_0$. In this case, we can find $\hat{g}$ such that $0 \leq \hat{g} \leq g$, $x_0 \notin {\rm supp}(\hat{g})$,
and $x_0 \in {\rm Int}(\hat{G}_t)$, where $ \hat{G}_t = \{x \in \mathbb{R}^n: d(x, {\rm supp}(\hat{g})) \leq t\}$. Then, as shown above, there exist $t_0 \in (0, t)$
and $\delta > 0$ such that
 $u^\varepsilon(t_0,x,y) \geq \hat{u}^\varepsilon(t_0,x,y) \geq 
\exp(-\varepsilon^{-1} \eta)$ for all sufficiently small~$\varepsilon$  and $\|x - x_0\|\leq \delta$, where $\hat{u}$ is the solution with the initial data $\hat{g}$. Thus we have proved that (\ref{ebel}) holds for some $t_0 \in (0, t)$.

 Consider now the diffusion process $(X^{\varepsilon}_t, Y^{\varepsilon}_t)$ starting at $(x, y)$ such that $\|x - x_0\| \leq \delta/2$.
Suppose $\eta' > 0$ is fixed. Let 
\[
\tau = \min(t - t_0, \inf\{s: u^\varepsilon(t - s,  X^{\varepsilon}_s, Y^{\varepsilon}_s) \geq 1 - \eta'\}).
\] 
Then
\[
u^\varepsilon(t,x,y) \geq
\]
\[
\mathrm{E}_{(x,y)} \left( u(\tau,X^{\varepsilon}_\tau, Y^{\varepsilon}_\tau) \exp(\varepsilon^{-1} \int_0^{\tau} c(Y^{\varepsilon}_s, 
u^\varepsilon(t-s, X^{\varepsilon}_s, Y^{\varepsilon}_s) ds),~\sup_{s \in [0,\tau]} \|X^{\varepsilon}_s - x\| \leq 
\frac{\delta}{2} \right).
\]
On the event $\{\tau = t - t_0,~ \sup_{s \in [0,\tau]} \|X^{\varepsilon}_s - x\| \leq 
\frac{\delta}{2}\}$, we have 
\[
u(\tau,X^{\varepsilon}_\tau, Y^{\varepsilon}_\tau) \exp(\varepsilon^{-1} \int_0^{\tau} c(Y^{\varepsilon}_s, 
u^\varepsilon(t-s, X^{\varepsilon}_s, Y^{\varepsilon}_s) ds) \geq \exp(-\varepsilon^{-1} \eta) \exp(\varepsilon^{-1} \hat{c}),
\]
where $\hat{c} = \min_{k = 1,2} \inf_{u \in [0, 1 -\eta']} (c^k(u))$. The right hand side can be made larger than one by
selecting a sufficiently small $\eta$.  On the event $\tau < t - t_0$, 
\[
u(\tau,X^{\varepsilon}_\tau, Y^{\varepsilon}_\tau) \exp(\varepsilon^{-1} \int_0^{\tau} c(Y^{\varepsilon}_s, 
u^\varepsilon(t-s, X^{\varepsilon}_s, Y^{\varepsilon}_s) ds) \geq 1 - \eta'.
\]
Finally,
\[
\lim_{\varepsilon \downarrow 0} \mathrm{P}_{(x,y)} ( \sup_{s \in [0,\tau]} \|X^{\varepsilon}_s - x\| \leq 
\frac{\delta}{2}) = 1.
\]
Therefore,
\[
u^\varepsilon(t,x,y)  \geq  1 - 2\eta'
\]
for all sufficiently small $\varepsilon$ and $(x, y) \in B_{\delta/2}(x_0) \times [0,1]$. Extracting a finite covering of $K$ by such domains, we see
that the estimate holds for $(x, y) \in K$ for for all sufficiently small $\varepsilon$. Since $\eta'$ was arbitrary, this implies the statement of the theorem.
\qed

\section{The case of $x$-dependent coefficients.} \label{lin1bb}

\subsection{Asymptotics of solutions to linear equations.}  \label{Asle}
In this section, we again consider the linear version of the Cauchy problem (\ref{rde1it}), but now allow the coefficients
$a^k,  \alpha^k, c^k$, $k =1,2$, to depend on~$x$. Thus $u^\varepsilon$ satisfies
\begin{equation} \label{rde1itlin22}
\frac{\partial u^\varepsilon}{\partial t} = \frac{\varepsilon}{2} \sum_{i,j =1}^n a_{ij}(x, y) \frac{\partial^2 u^\varepsilon}{\partial x_i \partial x_j}  + 
\frac{\varepsilon^{-\beta}}{2} \alpha(x, y)
\frac{\partial^2 u^\varepsilon}{\partial y^2} + \frac{1}{\varepsilon} c(x,y)u^\varepsilon,~~
t > 0,~x \in \mathbb{R}^n,~y \in (0,1) \setminus \{m\}. 
\end{equation}
\[
\frac{\partial u^\varepsilon}{\partial y}|_{y=0,1} = 0,~~~u^\varepsilon(0,x,y) = g(x).
\]
Recall that the pair of processes $(X^{\varepsilon}_t, Y^{\varepsilon}_t)$ starting at $(x,y)$ has been defined in Section~\ref{introd}.
In the case considered in Section~\ref{lint}, the main contribution to the expectation in (\ref{fkfo1}) comes from the event that the trajectoris 
of $X^{\varepsilon}_t$ (starting at $x$) closely follow the  linearly parametrized segment connecting $x$ with $x'$, where $x'$ is one of the points in $G_0$. Since the coefficients were spatially homogeneous, the contribution to
the expectation from such an event depended only on the difference between $x$ and $x'$. Now there is an optimal path $\varphi: [0,t] \rightarrow \mathbb{R}^n$
such that the trajectories of $X^{\varepsilon}_t$ following in its vicinity give the main contribution to the expectation. The shape of the
path depends on both initial point $x$ and the final point $x'$, and will be determined by examining the behavior of the slow component jointly with the distribution of the fast component (when tracking the fast component, all the points of $(0,m)$ can be identified, as well as all the points of $(m, 1)$, i.e., we can view
the fast component as a process with just two distinct values).

Let $K = \{1,2\}$. Consider the random
occupation measure on $( K  \times [0,t],  \mathcal{B}(K)  \times \mathcal{B}([0,t]))$:
\[
\nu^{\varepsilon}_{t,x,y} (\{1\} \times \Delta  ) = \int_\Delta \chi_{[0,m)} (Y^{\varepsilon}_s) ds,~~~
\nu^{\varepsilon}_{t,x,y} (\{2\} \times \Delta  )  = \int_\Delta \chi_{(m,1]} (Y^{\varepsilon}_s) ds,
\]
where $\Delta \in \mathcal{B}([0,t])$ and $(X^{\varepsilon}_0, Y^{\varepsilon}_0) = (x,y)$. 
The space of measures on $( K  \times [0,t],  \mathcal{B}(K)  \times \mathcal{B}([0,t]))$ whose marginals $\nu_s$, $s \in [0,t]$, are probability measures on $K$,
will be denoted by $\mathcal{M}$. It is endowed with the Levy-Prohorov distance denoted by $\rho$.  

Let $\mathcal{C}$ be the space of continuous functions on $[0,t]$ endowed with the distance $d$. Thus $X^{\varepsilon}$ can be viewed as a random
element of $\mathcal{C}$.

For $x, v \in \mathbb{R}^n$, define
\[
{R}(p, x, v) = \frac{1}{2} ((p_1 a^1(x) + p_2 a^2(x))^{-1} v, v).
\]

For $\varphi \in \mathcal{C}$ and $\nu \in \mathcal{M}$, define
\[
\bar{R}( \nu, \varphi) = \int_0^t {R}(\nu_s, {\varphi}_s, \dot{\varphi}_s) ds,
\]
For $f \in C([0,1])$, let $H^x(f)$ be the top eigenvalue 
of the operator $L^x_f u = \frac{1}{2} \alpha(x,y)
 u''(y) + fu$ (with the gluing condition at $y=m$ and reflection at the end points).  Let $\pi(x)$ be the invariant measure
for the process governed by this operator. 
For $\mu \in \mathcal{M}_{[0,1]}$, define 
\[
I^x(\mu) = \sup_{f \in C([0,1])} (\int_0^1 f d \mu - H^x(f)). 
\]
For $x \in \mathbb{R}^n$ and $p \in J$ , define 
\[
S(p,x) = \inf_{\mu: p_\mu = (p_1, p_2)} I^x(\mu).
\]
For $\varphi \in \mathcal{C}$ and $\nu \in \mathcal{M}$, define
\[
\bar{S}(\nu,\varphi) = \int_0^t S(\nu_s, \varphi_s) ds.
\]
Let $\tilde{\Lambda}^\varepsilon_{t,x,y}$, be
the measure on $(\mathcal{M}  \times \mathcal{C}, \rho \times d )$ induced by $(\nu^\varepsilon_{t,x,y}, X^\varepsilon)$ (with $X^\varepsilon_0 = x$ and $Y^\varepsilon_0 = y$). Note that $\tilde{\Lambda}^\varepsilon_{t,x,y}$ also depends on $\beta$ because of the dependence of $Y^\varepsilon_t$ on $\beta$. 
 The following
theorem is proved, in a somewhat different form, in~\cite{Lip}.

\begin{theorem} \label{ldp} If $\beta = 1$, 
the family $\tilde{\Lambda}^\varepsilon_{t,x,y}$ obeys the large deviations principle with the action functional 
\[
L(\nu , \varphi) =  \bar{R}(\nu, \varphi) + \bar{S}(\nu, \varphi),
\]
uniformly in $(x, y)$ on every compact $K \subset \mathbb{R}^n \times [0,1]$.

 If $\beta > 1$, the family $\tilde{\Lambda}^\varepsilon_{t,x,y}$ obeys the large deviations principle with the action functional 
\[
L(\nu , \varphi) =  \begin{cases} \bar{R}(\tilde{\nu}^\varphi, \varphi) ~~~ if~ \nu = \tilde{\nu}^\varphi, \\ \infty ~~~~~~~~~ otherwise, \end{cases} 
\]
uniformly in $(x, y)$ on every compact $K \subset \mathbb{R}^n \times [0,1]$, where $\tilde{\nu}^\varphi$ is such that $\tilde{\nu}^\varphi_s = \pi(\varphi_s)$ for each $s$. 

If $-1 < \beta < 1$, the family $\tilde{\Lambda}^\varepsilon_{t,x,y}$ obeys the large deviations principle with the action functional
\[
L(\nu , \varphi) =  \inf_{\nu' \in \mathcal{M}}  \bar{R}(\nu', \varphi),
\]
uniformly in $(x, y)$ on every compact $K \subset \mathbb{R}^n \times [0,1]$.
%
%
%
\end{theorem}

Define
\[
{T}(p, x) =  p_1 c^1(x) + p_2 c^2(x)
\]
and
\[
\bar{T}(\nu, \varphi) = \int_0^t T(\nu_s, \varphi_s ) ds.
\]
Let $\mathcal{C}(x, x') = \{\varphi \in \mathcal{C}:  \varphi(0) = x, \varphi(t)  =x' \}$. 
For $\beta = 1$, define
\begin{equation} \label{mo11} 
\lambda(t,x, x') =  \sup_{\varphi \in \mathcal{C}(x, x')} \sup_{\nu \in \mathcal{M}} (\bar{T}(\nu, \varphi) - \bar{S}(\nu, \varphi) - \bar{R}(\nu, \varphi)). 
\end{equation}
The expression for $\lambda(t,x, x')$ in the case when $\beta > 1$ is
\begin{equation} \label{ltmo22}
\lambda(t,x, x') = \sup_{\varphi \in \mathcal{C}(x, x')} (\bar{T}(\tilde{\nu}^\varphi, \varphi)  - \bar{R}(\tilde{\nu}^\varphi, \varphi)).
\end{equation} 
Finally, in the case when $-1 < \beta < 1$, define
\begin{equation} \label{gtmo33}
\lambda(t,x, x') =  \sup_{\varphi \in \mathcal{C}(x, x')} \sup_{\nu \in \mathcal{M}} (\bar{T}(\nu, \varphi)  - \bar{R}(\nu, \varphi)). 
\end{equation}
\begin{theorem} \label{linle22}
Under the above assumptions, 
\begin{equation} \label{limyy22}
\lim_{\varepsilon \downarrow 0} \varepsilon \ln ({u}^\varepsilon(t,x,y)) =  \sup_{x' \in G_0} \lambda(t, x, x')
\end{equation}
uniformly on every compact $K \subset (0,\infty) \times \mathbb{R}^n \times [0,1]$, where $\lambda$ is given by (\ref{mo11}) if $\beta = 1$, by 
(\ref{ltmo22}) if $\beta > 1$, and (\ref{gtmo33}) if $-1 < \beta < 1$. 
\end{theorem}
\proof The proof is similar to that of Theorem \ref{linle22}. The main difference is that
in (\ref{fkfo1}) we were able to represent ${u}^\varepsilon(t,x,y)$ is terms of a repeated integral with respect to the measures
$M^{x, p, \varepsilon}_t$ and ${\Lambda}^{\varepsilon}_{t,y}$. Now, we'll instead use the measure  $\tilde{\Lambda}^\varepsilon_{t,x,y}$
 on $(\mathcal{M}  \times \mathcal{C}, \rho \times d )$.
Fix $t > 0$.  By the Feynman-Kac formula, 
\[
u^\varepsilon(t,x,y) = \mathrm{E}_{(x,y)} \left( g(X^\varepsilon_t) \exp(\varepsilon^{-1} \int_0^t c(X^\varepsilon_s, Y^\varepsilon_s) ds) \right) =
\]
\[
\int_{\mathcal{M}  \times \mathcal{C}} g(\varphi_t) \exp\left(\varepsilon^{-1} \int_0^t \sum_{k\in K} c^k(\varphi(s)) \nu_s(k)d s\right) 
d \tilde{\Lambda}^\varepsilon_{t,x,y}(\nu, \varphi).
\]  
By (\ref{actf}) and Theorem~\ref{ldp}, 
\[
\lim_{\varepsilon \downarrow 0} \varepsilon \ln ({u}^\varepsilon(t,x,y))  = \sup_{x' \in G_0}  
 \sup_{\varphi \in \mathcal{C}(x, x')} \sup_{\nu \in \mathcal{M}} (\bar{T}(\nu, \varphi) - L(\nu , \varphi)).
\]
In each of the cases, $\beta = 1$, $\beta > 1$, and $-1 < \beta < 1$, we can insert the expression for $L(\nu , \varphi)$ provided in Theorem~\ref{ldp} into
the right hand side of the last formula. Thus we obtain that (\ref{limyy22}) holds uniformly in $(x, y) \in \bar{K} \times [0,1]$, where  $\bar{K} \subset
\mathbb{R}^n$ is compact. The  uniform convergence on $K \subset (0,\infty) \times \mathbb{R}^n \times [0,1]$ can be justified in the same as in Theorem~\ref{linle22}.
\qed

\subsection{Asymptotics of solutions to reaction-diffusion equations.} 

As in Section~\ref{rde1}, here  we consider the Cauchy problem for the reaction-diffusion equation~(\ref{rde1it}), but now we allow $a^k, \alpha^k, c^k$, $k =1,2$, to
 depend on~$x$. 
Let $\tilde{c}^1(x) = c^1(x,0)$, $\tilde{c}^2(x) = c^2(x,0)$. 
Consider the linear problem (\ref{rde1itlin22})
with $c^1$, $c^2$ replaced by $\tilde{c}^1$, $\tilde{c}^2$. Let $\lambda(t,x,x')$ be
given by (\ref{mo11}) if $\beta = 1$, by 
(\ref{ltmo22}) if $\beta > 1$, and (\ref{gtmo33}) if $-1 < \beta < 1$. 
Define
\[
G_t = \{x \in \mathbb{R}^n: \lambda(s, x, \varphi(s)) \geq 0~{\rm for}~{\rm all}~s \in [0,t],~{\rm for}~{\rm some}~ \varphi \in \mathcal{C}~{\rm with}~ 
\varphi(0) = x,~\varphi(t) \in G_0\}. 
\]
This set (or, rather, $G_t \times [0,1]$) is the multi-layer analogue of the set $\{x: {V}(t,x) = 0 \}$, where $V$ was defined in the Introduction.
It is not difficult to show that if $\tilde{c}^1(x) \equiv \tilde{c}^2(x)  \equiv {\rm const} $, then the growth of $G_t$ obeys the Huygens principle 
with respect to a certain non-homogeneous metric. The metric satisfies $d(x, x')  = \inf \{t \geq 0: \lambda(t, x, x') \geq 0\}$, where
$\lambda$ was defined in Section~\ref{Asle}. 
\begin{theorem} \label{th44} If $u^\varepsilon(t,x,y)$ is the solution of (\ref{rde1it}) and $c$ is of KPP type, then, for each $t > 0$, 
\[
\lim_{\varepsilon \downarrow 0} u^\varepsilon(t,x,y)) = 0
\]
uniformly on every compact $K \subset (\mathbb{R}^n \setminus G_t) \times [0,1]$, and
\[
\lim_{\varepsilon \downarrow 0} u^\varepsilon(t,x,y)) = 1
\]
uniformly on every compact $K \subset {\rm Int}(G_t) \times [0,1]$.
\end{theorem}
Before we proceed with the proof of this theorem, let us discuss an example. Let $\beta = 1$. 
Assume that $n = 1$ and $G_0 = [-2,-1]$. Suppose that $a^1 = a^2 \equiv 1$. Let us take $\alpha^1(x) = \alpha^2(x) = \delta^{-1}$ for
$x < -\delta$, $\alpha^1(x) = \alpha^2(x) = \delta$ for $x > \delta$, and $\alpha^1(x) = \alpha^2(x) \in [\delta, \delta^{-1}]$ for $x \in [-\delta, \delta]$. 
Assume that $\tilde{c}^1 (x) \equiv \delta$, while $\tilde{c}^2(x) \equiv 1$.
We also assume that $m = 2/3$, i.e., the first layer is twice as thick as the second one. Optimizing over
the time $s \in [0,t]$ that a trajectory $\varphi$ spends to the right of the origin, from (\ref{mo11}) we obtain that
\[
\lim_{\delta \downarrow 0} \lambda(t, 0, -1) = \sup_{s \in [0,1]} ( s + \frac{t-s}{3} - \frac{1}{2(t-s)}).
\]
The supremum in the right hand side is achieved when $s = t - \sqrt{3}/2$. The right hand side is positive if and only if $t > 2/\sqrt{3}$. Moreover,
if $t > 2/\sqrt{3}$ and $\delta > 0$ is sufficiently small then the trajectory $\varphi$ such that $\varphi(\tau) = \delta$ for $0 \leq \tau \leq t - \sqrt{3}/2$,
$\varphi(\tau) = \delta - \frac{2(1 + \delta)}{\sqrt{3}} (\tau - (t-\frac{\sqrt{3}}{2}))$  for $\tau \in [t - \sqrt{3}/2,t]$, has the property that
$\lambda(\tau, \delta, \varphi(\tau)) \geq 0$  for  all $\tau \in [0,t]$. Thus, if $t > 2/\sqrt{3}$, then $\delta \in G_t$ for all sufficiently small $\delta$. 
At the same time, it is not difficult to check that
\[
\lim_{\delta \downarrow 0} \lambda(2/\sqrt{3}, z, -1) < 0
\]
if $z < 0$ and $|z|$ is sufficiently small. Thus, $z \notin G_t$ if $t > 2/\sqrt{3}$ and $t - 2/\sqrt{3}$ and $\delta$ are sufficiently small. This demonstrates
that the interface jumps at some time prior to $t$.
\\
\\
{\it Proof of Theorem~\ref{th44}.} The main difference from the proof of Theorem~\ref{nnln} is that now there may exist $x \notin G_t$ such that $\lambda(t,x,x') > 0$
for some $x' \in G_0$. This is due to the fact that, in general, $\lambda(t,x,x') > 0$ does not imply that 
there is $\varphi \in \mathcal{C}(x, x')$ such that $\lambda(s,x,\varphi(s)) > 0$ for each $s \in [0,t)$. Thus, comparison with the solution of the corresponding linear
equation is not immediately available to establish the first statement of the theorem. However, once we prove the first statement, the proof that 
$\lim_{\varepsilon \downarrow 0} u^\varepsilon(t,x,y) = 1$ uniformly on every compact inside ${\rm Int}(G_t) \times [0,1]$ is similar to that in Theorem~\ref{nnln},
and thus we focus on proving that $\lim_{\varepsilon \downarrow 0} u^\varepsilon(t,x,y) = 0$
uniformly on  $K \subset (\mathbb{R}^n \setminus G_t) \times [0,1]$.

Let $(x_0, y_0) \in K$. Given $\delta > 0$, let $U_\delta \subseteq \mathbb{R}^n$ be the
set defined as follows: $x' \in U_\delta$ if there exists  $\varphi \in \mathcal{C}$ such that $\varphi(0) = x_0$, $\varphi(t) = x'$, and
$ \lambda(s + \delta, x_0, \varphi(s)) \geq 0 $ for $0 \leq s \leq t$.  By the definition of $G_t$, there is $\delta > 0$ such that $U_\delta $ and $ G_0$ are disjoint. 
Observe that $x_0 \in U_\delta$ and $\lambda(t+\delta, x_0, x') = 0$ for $x' \in \partial U_\delta$. Let $\tilde{U}_\delta$ be a bounded
 domain with a smooth boundary such that
$x_0 \in \tilde{U}_\delta \subset U_\delta$ and $\lambda(t+\delta/2,x_0, x') \leq 0$ for $x' \in \partial \tilde{U}_\delta$. Moreover, we can choose this domain in such 
a way that, for some $r > 0$ and all $x$ satisfying $\|x - x_0\| \leq r$,  $x \in \tilde{U}_\delta \subset U_\delta$ and $\lambda(t+\delta/2,x, x') \leq 0$ for $x' \in \partial \tilde{U}_\delta$.

Let $\tau = \inf\{t: X^\varepsilon_t \in \partial \tilde{U}_\delta\}$. 
Let $\tilde{u}^\varepsilon$ solve the following initial-boundary value problem:
\begin{equation} \label{rde1itlin22ib}
\frac{\partial \tilde{u}^\varepsilon}{\partial t} = \frac{\varepsilon}{2} \sum_{i,j =1}^n a_{ij}(x, y) \frac{\partial^2 
\tilde{u}^\varepsilon}{\partial x_i \partial x_j}  + 
\frac{\varepsilon^{-\beta}}{2} \alpha(x,y)
\frac{\partial^2 \tilde{u}^\varepsilon}{\partial y^2} + \frac{1}{\varepsilon} \tilde{c}(x,y)\tilde{u}^\varepsilon,~~
t > 0,~x \in \tilde{U}_\delta,~y \in (0,1) \setminus \{m\}. 
\end{equation}
\[
\frac{\partial \tilde{u}^\varepsilon}{\partial y}|_{y=0,1} = 0,~~~\tilde{u}^\varepsilon(0,x,y) = 0,~~x \in \tilde{U}_\delta;~~~\tilde{u}^\varepsilon(t,x,y) = 1,~~
x \in \partial \tilde{U}_\delta,
\]
 where the solution is defined using the Feynman-Kac formula, 
\[
\tilde{u}^\varepsilon(t,x,y) = \mathrm{E}_{(x,y)} \left( \chi_{\tau \leq t}  \exp(\varepsilon^{-1} \int_0^\tau c(X^\varepsilon_s, Y^\varepsilon_s) ds) \right).
\]
It is clear that $u^\varepsilon(t,x,y) \leq \tilde{u}^\varepsilon(t,x,y)$. 
The desired statement follows from the fact that $\lim_{\varepsilon \downarrow 0} \tilde{u}^\varepsilon(t,x,y) = 0$ uniformly for $(x,y) \in B_r(x_0) \times [0,1]$. 
To show the latter, we modify the proof of Theorem~\ref{linle22}.
For $\varphi \in \mathcal{C}$, let $\tau(\varphi) = \inf\{t: \varphi(t) \in \partial \tilde{U}_\delta\}$. Let
\[
\hat{T}(\nu, \varphi) =  \begin{cases} \int_0^{\tau(\varphi)} T(\nu_s, \varphi_s ) 
ds ~~~~~ {\rm if}~ \tau(\varphi) \leq t \\ 0 ~~~~~~~~~~~~~~~~~~~~~~~~~{\rm otherwise}. \end{cases} 
\]
Then 
\[
\tilde{u}^\varepsilon(t,x,y) = \mathrm{E}_{(x,y)} \left( \chi_{\tau \leq t}  \exp(\varepsilon^{-1} \int_0^\tau c(X^\varepsilon_s, Y^\varepsilon_s) ds) \right) =
\]
\[
\int_{\mathcal{M}  \times \mathcal{C}} \chi_{\tau(\varphi) \leq t} \exp\left(\varepsilon^{-1} \int_0^{\tau(\varphi)} \sum_{k\in K} c^k(\varphi(s)) \nu_s(k)d s\right) 
d \tilde{\Lambda}^\varepsilon_{t,x,y}(\nu, \varphi).
\]  
By (\ref{actf2}), for each $\eta > 0$, 
\[
\limsup_{\varepsilon \downarrow 0} \varepsilon \ln (\tilde{u}^\varepsilon(t,x,y))  \leq 
 \sup_{\varphi \in \mathcal{C}, \varphi(0) = x} \sup_{\nu \in \mathcal{M}} (\hat{T}(\nu, \varphi) - L^\eta(\nu , \varphi)),
\]
where $L^\eta(\nu , \varphi) = \inf_{\rho(\nu',\nu)+d(\varphi',\varphi) \leq \eta} L^\eta(\nu' , \varphi')$. Pick $\bar{\nu}$ and $\bar{\varphi}$
such that 
\[
\hat{T}(\bar{\nu}, \bar{\varphi}) - L^\eta(\bar{\nu}, \bar{\varphi}) \geq \sup_{\varphi \in \mathcal{C}, \varphi(0) = x}\sup_{\nu \in \mathcal{M}} (\hat{T}(\nu, \varphi) - L^\eta(\nu , \varphi)) - \delta',
\]
where $\delta' > 0$ will be selected below. Then there are $\nu'$ and $\varphi'$ with $\rho(\nu',\bar{\nu})+d(\varphi',\bar{\varphi}) \leq \eta$ such that
\[
\hat{T}(\bar{\nu}, \bar{\varphi}) - L(\nu', \varphi') \geq \sup_{\varphi \in \mathcal{C}, \varphi(0) = x}\sup_{\nu \in \mathcal{M}} (\hat{T}(\nu, \varphi) - L^\eta(\nu , \varphi)) - 2\delta'.
\]
For $s \in [0,t]$, let $L_s(\nu, \varphi)$ be defined as $L(\nu, \varphi)$ on the interval $[0,s]$. Observe that $L_s(\nu, \varphi) \leq L(\nu,\varphi)$. Thus
\[
\hat{T}(\bar{\nu}, \bar{\varphi}) - L_{\tau(\bar{\varphi})}(\nu', \varphi') \geq \sup_{\varphi \in \mathcal{C}, \varphi(0) = x}\sup_{\nu \in \mathcal{M}} (\hat{T}(\nu, \varphi) - L^\eta(\nu , \varphi)) - 2\delta'.
\]
If $\eta$ is sufficiently small, this implies that
\[
\int_0^{\tau(\bar{\varphi})} T(\nu'_s, \varphi'_s ) ds  - L_{\tau(\bar{\varphi})}(\nu', \varphi') \geq \sup_{\varphi \in \mathcal{C}, \varphi(0) = x} \sup_{\nu \in \mathcal{M}} (\hat{T}(\nu, \varphi) - L^\eta(\nu , \varphi)) - 3\delta'.
\]
Since $\lambda(t+\delta/2,x, x') \leq 0$ for $x' \in \partial \tilde{U}_\delta$, the left hand side can be made smaller than $-4 \delta'$ by choosing a 
sufficiently small $\delta'$ and $\eta$, which shows that $\lim_{\varepsilon \downarrow 0} \tilde{u}^\varepsilon(t,x,y) = 0$ uniformly for $(x,y) \in B_r(x_0) \times [0,1]$. 
\qed

\section{Remarks and generalizations.} \label{remgen} 
{\bf 1.} Consider a composite consisting of periodic inclusions in homogeneous media (see  Figure 2). Suppose, for brevity, that the system is invariant with respect
to shifts of size one in each variable. Assume that the inclusions are domains with diameter $\delta$, and that each inclusion contains
a ball of diameter $\delta/N$, where $N > 1$ is constant. If the nonlinear term $f(x,u) = u c(x, u)$ is of KPP type, the growth of the domain where
$u^\varepsilon(t,x)$ is close to one for $t \gg  1$, with fixed  $\varepsilon$ and $\delta$,  has been described in \cite{F3}, Ch. 7. 

Let $u^{\varepsilon,\delta}(t,x)$ be the solution of problem (\ref{rde1itone}) in this medium. Equation (\ref{eqw3}), together with large deviations
estimates, allows one to describe the limiting behavior of $u^{\varepsilon,\delta(\varepsilon)}(t,x)$ as $\varepsilon \ll \delta(\varepsilon) \ll 1$. Suppose, for
brevity, that $n = 1$, $a_{11}(x) \equiv 1$, $c(x,0) = c_0$ outside the inclusions, and $c(x,0) = c_1$ inside the inclusions, where $c_0$ and $c_1$ are constants. 
Moreover, assume that $c_1 > 2 c_0$, $G_0 = (-\infty, 0]$, and that the interval $[-\delta, 0]$ coinsides with  one of the inclusions. 
Then, according to \cite{F2} (see
also \cite{F3}, Ch 6.2),  as $0 < \varepsilon \ll 1$, the interface first moves to the right with the speed $\sqrt{2 c_0}$, and then jumps to $x_1  = 
1-\delta$  at the time $T_0 = (1- \delta) \sqrt{2(c_1 - c_0)}/ c_1 < 1/ \sqrt{2 c_0}$. 

This implies that the average speed of the expansion of the region where $u^{\varepsilon,\delta}$ is close to one will be arbitrarily large if $c_1$ is large enough.
Moreover, the choice of $c_1$ providing the rapid expansion is independent of the fraction $\delta$ of the inclusions in the composite. This effect, under some
additional assumptions, is preserved for other types of nonlinearities. For instance, if $f(x,u)$ is a bistable nonlinearity outside of the inclusions (such as 
$f(x,u) = u(1-u)(u-\lambda)$ with $\lambda \in (0,1)$ and $\int_0^1 f(x,u) du > 0$),  then one has this acceleration of the expansion due to 
the  fact that the  nonlinearity in the inclusions is of KPP type. Moreover, this effect is preserved for non-periodic inclusions if, in a sense, their  sizes and their fraction satisfy certain bounds from below.    

{\bf 2.} Asymptotic problems for RDEs with other types of nonlinearities can be also considered in the layered composites. For example, if $f(x,u)$ is a bistable nonlinear term (as above) and $\beta > 1$ in equation (\ref{rde1it}), the propagation of the region where $u^{\varepsilon,\delta}$ is close to one can be
described by the Huygens principle in the average Reimannian metric as considered in Section~\ref{lin1bb}. The interface motion in  the  bistable case
always has a local nature, and the constant velocity in   this metric  is defined as the speed of the front in a one-dimensional space-homogeneous medium 
(compare with \cite{Gart1}). The proof of this statement can be derived from the bounds obtained in \cite{Gart1}. 

{\bf 3.} One can consider RDEs where the reaction occurs just on the surface $\{y = m\}$ dividing the layers. In this case, problem (\ref{rde1it}) should be
modified: the nonlinear term should be excluded from the equation, and the gluing condition has the form
\[
\frac{\partial^+ u^\varepsilon(t,x,y)}{\partial y} |_{y = m} - \frac{\partial^- u^\varepsilon(t,x,y)}{\partial y} |_{y = m} = 
\frac{1}{\varepsilon} c(x, u^\varepsilon) u^\varepsilon,
\]
where the differentials $\partial^+$ and $\partial^-$ mean that the derivatives are calculated when $y$ approaches $m$ from above and below, respectively. The modified
Feynman-Kac formula  in this case gives the following equation for $u^\varepsilon(t,x,y)$:
\[
u^\varepsilon(t,x,y) =  \mathrm{E}_{(x,y)} \left( g(X^{\varepsilon}_{t}) \exp(\varepsilon^{-1} \int_0^{t} c(X^\varepsilon_s, Y^{\varepsilon}_s, 
u^\varepsilon(t -s, X^{\varepsilon}_s, Y^{\varepsilon}_s) ) dL^\varepsilon_s )\right),
\]
where $L^\varepsilon_t$ is the local time of the process $(X^\varepsilon_t, Y^{\varepsilon}_t)$ on the surface $\{y = m\}$. Here, additional diffuculties arise 
due to the large deviations for the local time. If the diffusion coefficients are continuous on the surface $\{y = m\}$, the problem can be studied similarly to 
\cite{Fs}. (Now, however, the action functional is more sophisticated than the one considered in \cite{Fs}.) \be

{\bf 4.}  Finally, we would like to mention that effects caused by random thickness of the layers, random distribution of inclusions, as well as other types of 
underlying stochastic transport, can also be studied using large deviation asymptotics. We will address some of these problems in a different paper. 
\\
\\

\noindent {\bf \large Acknowledgements}:  While working on this
article, M. Freidlin was supported by NSF grant DMS-1411866
and L. Koralov was supported by ARO grant W911NF1710419.
\\
\\

\end{document}